\documentclass[11pt]{article}
\usepackage{amscd,amsmath,amssymb}
\usepackage[mathscr]{eucal}

\makeatletter \@addtoreset{equation}{section}\makeatother

\title{\bf QUANTUM MATRIX BALL: THE CAUCHI-SZEG\"O KERNEL AND THE SHILOV
BOUNDARY}

\author{L. Vaksman} \date{}

\newtheorem{theorem}{Theorem}[section]
\newtheorem{lemma}[theorem]{Lemma}
\newtheorem{proposition}[theorem]{Proposition}

\begin{document}
\large
\maketitle

\makeatletter
\let\@thefnmark\relax
\@footnotetext{This research was supported in part by Award No UM1-2091 of
the US Civilian Research \& Development Foundation}
\makeatother

\begin{abstract}

This work produces a q-analogue of the Cauchi-Szeg\"o integral
representation that retrieves a holomorphic function in the matrix ball from
its values on the Shilov boundary. Besides that, the Shilov boundary of the
quantum matrix ball is described and the $U_q
\mathfrak{su}_{m,n}$-covariance of the $U_q \mathfrak{s}(\mathfrak{u}_m
\times \mathfrak{u}_n)$-invariant integral on this boundary is established.
The latter result allows one to obtain a q-analogue for the principal
degenerate series of unitary representations  related to the Shilov
boundary of the matrix ball.
\end{abstract}

\section{Introduction}

Bounded symmetric domains form a favorite subject of research in geometry,
function theory of several complex variable, and non-commutative harmonic
analysis. The point is that they are simplest among non-compact homogeneous
spaces of real semi-simple Lie groups.

Quantum analogues (q-analogues) of bounded symmetric domains were introduced
in \cite{SV} via replacement of the ordinary Lie groups with their quantum
analogues \cite{CP}. A simple example of classical bounded symmetric domain
is the unit ball $\mathbb{U}=\{\mathbf{z}\in
\mathrm{Mat}_n|\,\mathbf{zz}^*<I \}$ in the space of complex $n \times n$
matrices. The present work studies a q-analogue of this matrix ball. The
notions of Shilov boundary and Cauchi-Szeg\"o kernel are associated to this
subject. Just as in the case $q=1$, a holomorphic function in the quantum
matrix ball is determined unambiguously by its restriction onto the Shilov
boundary and could be retrieved via the Cauchi-Szeg\"o integral
representation.

The proof of the latter result involves some properties of an invariant
integral on the Shilov boundary. We also use these properties to produce the
principal degenerate series of unitary representations of the quantum group
$SU_{n,n}$. The appendix contains a discussion of generalizations of the
main result onto the case of rectangular matrices.

In what follows we assume $\mathbb{C}$ to be the ground field, all the
algebras are unital, and $q \in(0,1)$.

The author is grateful to S. Sinel'shchikov for pointing out a gap in a
previous version of the proof of proposition \ref{ng} and an assistance with
translating this text into English.

\bigskip

\section{A construction of the Shilov boundary}

Among the best known 'quantum' algebras one should mention the q-analogue
$\mathbb{C}[\mathrm{Mat}_n]_q$ of the algebra of holomorphic polynomials on
the space of matrices $\mathrm{Mat}_n$. This algebra is given by its
generators $z_a^\alpha$, $a,\alpha=1,2,\ldots,n$, and the following
relations:
\begin{align}
\label{zaa1} z_a^\alpha z_b^\beta-qz_b^\beta z_a^\alpha &=0, & a=b &
\;\&\;\alpha<\beta & {\rm or}\qquad a<b \;\&\; \alpha=\beta \\ \label{zaa2}
z_a^\alpha z_b^\beta-z_b^\beta z_a^\alpha &=0,& \alpha<\beta & \;\&\;a>b &
\\ \label{zaa3} z_a^\alpha z_b^\beta-z_b^\beta z_a^\alpha
&=(q-q^{-1})z_a^\beta z_b^\alpha, & \alpha<\beta & \;\&\; a<b. &
\end{align}

The work \cite{SSV1} presents a definition of the $*$-algebra
$\mathrm{Pol}(\mathrm{Mat}_n)_q$, which is a q-analogue of the polynomial
algebra on the vector space $\mathrm{Mat}_n$. Its generators are
$z_a^\alpha$, $(z_a^\alpha)^*$, $a,\alpha=1,2,\ldots,n$, and the list of
relations consists of (\ref{zaa1}), (\ref{zaa2}), (\ref{zaa3}), together
with the commutation relations
\begin{equation}\label{rmr}
(z_b^\beta)^*z_a^\alpha=q^2
\sum_{a',b',\alpha',\beta'=1}^nR_{ba}^{b'a'}R_{\beta\alpha}^{\beta'\alpha'}
z_{a'}^{\alpha'}(z_{b'}^{\beta'})^*+(1-q^2)\delta_{ab}\delta^{\alpha\beta},
\end{equation}
where $\delta_{ab}$, $\delta^{\alpha \beta}$ being the Kronecker symbols and
$$
R_{ij}^{kl}=
\begin{cases}
q^{-1}, & i \ne j\;\&\;i=k\;\&\;j=l,\\ 1,& i=j=k=l,\\ -(q^{-2}-1),& i=j
\;\&\;k=l \;\&\; l>j,\\ 0 & \text{otherwise}.
\end{cases}
$$

It is well known that the Shilov boundary of the matrix ball $\mathbb{U}$ is
just the set $S(\mathbb{U})$ of all unitary matrices. Our intention is to
produce a q-analogue of the Shilov boundary for the quantum matrix ball.
Introduce the notation for the quantum minors of the matrix
$\mathbf{z}=(z_a^\alpha)$:
$$
(z^{\wedge k})_{\{a_1,a_2,\ldots,a_k
\}}^{\{\alpha_1,\alpha_2,\ldots,\alpha_k
\}}\stackrel{\mathrm{def}}{=}\sum_{s \in
S_k}(-q)^{l(s)}z_{a_1}^{\alpha_{s(1)}}z_{a_2}^{\alpha_{s(2)}}\cdots
z_{a_k}^{\alpha_{s(k)}},
$$
with $\alpha_1<\alpha_2<\ldots<\alpha_k$ and $l(s)$, $a_1<a_2<\ldots<a_k$,
and $l(s)$ being the number of inversion in $s \in S_k$.

It is well known that the quantum determinant
$$
\det \nolimits_q \mathbf{z}=(z^{\wedge n})_{\{1,2,\ldots,n
\}}^{\{1,2,\ldots,n \}}
$$
is in the center of $\mathbb{C}[\mathrm{Mat}_n]_q$. The localization of
$\mathbb{C}[\mathrm{Mat}_n]_q$ with respect to the multiplicative system
$(\det_q \mathbf{z})^\mathbb{N}$ is called the algebra of regular functions
on the quantum $GL_n$ and is denoted by $\mathbb{C}[GL_n]_q$.

\medskip

\begin{lemma}\label{inv}
i) There exists a unique involution $*$ in $\mathbb{C}[GL_n]_q$ such that
$$
(z_a^\alpha)^*=(-q)^{a+\alpha-2n}(\det \nolimits_q \mathbf{z})^{-1}\det
\nolimits_q \mathbf{z}_a^\alpha,
$$
with $\mathbf{z}_a^\alpha$ being the matrix derived from $\mathbf{z}$ via
deleting the line $\alpha$ and the column $a$.

ii) $(\det_q \mathbf{z})(\det_q \mathbf{z})^*=(\det_q \mathbf{z})^*(\det_q
\mathbf{z})=q^{-n(n-1)}$.
\end{lemma}

\smallskip

{\bf Proof.} The uniqueness of the involution $*$ is obvious. To prove the
existence, consider the $*$-algebra
$\mathbb{C}[U_n]_q=(\mathbb{C}[GL_n]_q,\star)$ of regular functions on the
quantum $U_n$ (see \cite{Koe}) and the automorphism $i:\mathbb{C}[GL_n]_q
\to \mathbb{C}[GL_n]_q$ given by $i:z_a^\alpha \mapsto
q^{\alpha-n}z_a^\alpha$, $a,\alpha=1,2,\ldots,n$. Obviously, $i^{-1}\star i$
is an involution. What remains is to demonstrate that
$$
i^{-1}\star i: z_a^\alpha \mapsto (-q)^{a+\alpha-2n}(\det \nolimits_q
\mathbf{z})^{-1}\det \nolimits_q \mathbf{z}_a^\alpha.
$$
This follows from the definition of $\star$ \cite{Koe}:
$$
(z_a^\alpha)^\star=(-q)^{a-\alpha}(\det \nolimits_q \mathbf{z})^{-1}\det
\nolimits_q \mathbf{z}_a^\alpha.
$$

The statement ii) of the lemma follows from a similar statement for the
involution $\star$. \hfill $\square$

\medskip

The $*$-algebra $\mathrm{Pol}(S(\mathbb{U}))_q=(\mathbb{C}[GL_n]_q,*)$ is a
q-analogue of the polynomial algebra on the Shilov boundary of the matrix
ball $\mathbb{U}$, as one can see from

\begin{theorem}\label{psi}
There exists a unique homomorphism of $*$-algebras
$\psi:\mathrm{Pol}(\mathrm{Mat}_n)_q \to \mathrm{Pol}(S(\mathbb{U}))_q$ such
that $\psi:z_a^\alpha \mapsto z_a^\alpha$, $a,\alpha=1,2,\ldots,n$.
\end{theorem}

\smallskip

We premise the proof of the theorem with two remarks. Firstly, the
homomorphism $\psi$ is a q-analogue of the operator which restricts the
polynomial onto the Shilov boundary. Secondly, we use in this work a purely
algebraic approach to producing the Shilov boundary; we do not try to
compare it to the analytic approach of the well known work by W. Arveson
\cite{Ar}.

\smallskip

{\bf Proof.} The uniqueness of the $*$-homomorphism $\psi$ is obvious.
Turn to the proof of its existence. To produce $\psi$, we need auxiliary
$*$-algebras $\mathrm{Pol}(\mathrm{Pl}_{n,2n})_{q,x}$ and $\mathscr{F}_x$,
together with the embeddings of algebras $\mathrm{Pol}(\mathrm{Mat}_n)_q
\hookrightarrow \mathrm{Pol}(\mathrm{Pl}_{n,2n})_{q,x}$,
$\mathrm{Pol}(S(\mathbb{U}))_q \hookrightarrow \mathscr{F}_x$. The crucial
point in the proof is a construction of a homomorphism of $*$-algebras
$\mathrm{Pol}(\mathrm{Pl}_{n,2n})_{q,x} \to \mathscr{F}_x$ which leads to
the commutative diagram
$$
\begin{CD}
\mathrm{Pol}(\mathrm{Pl}_{n,2n})_{q,x}@>>>\mathscr{F}_x \\ @AAA @AAA \\
\mathrm{Pol}(\mathrm{Mat}_n)_q @>>>\mathrm{Pol}(S(\mathbb{U}))_q.
\end{CD}
$$
To begin with, introduce the $*$-algebra
$\mathrm{Pol}(\mathrm{Pl}_{n,2n})_{q,x}$ 'of polynomials on the quantum
Pl\"ucker manifold'.

Let $\mathbb{C}[\mathrm{Mat}_{2n}]_q$ be the $*$-algebra of functions on the
quantum space of matrices $\mathrm{Mat}_{2n}$ determined by its generators
$\{t_{ij} \}_{i,j=1,2,\ldots,2n}$ and commutation relations similar to those
listed in (\ref{zaa1}) -- (\ref{zaa3}). Introduce the quantum minors
$$
t_{IJ}^{\wedge n}=\sum_{s \in
S_n}(-q)^{l(s)}t_{i_1j_{s(1)}}t_{i_2j_{s(2)}}\cdots t_{i_nj_{s(n)}}
$$
determined by pairs of $n$-element subsets of the form
$$I=\{(i_1,i_2,\ldots,i_n)|\:1 \le i_1<\ldots<i_n \le 2n \},$$
$$J=\{(j_1,j_2,\ldots,j_n)|\:1 \le j_1<\ldots<j_n \le 2n \}.$$

Consider the subalgebra in $\mathbb{C}[\mathrm{Mat}_{2n}]_q$ generated by
the elements $t_{\{1,2,\ldots,n \}J}^{\wedge n}$, $t_{\{n+1,n+2,\ldots,2n
\}J}^{\wedge n}$, with $\mathrm{card}(J)=n$. It is easy to present a full
list of relations between these generators. We follow \cite{SSV2} in
equipping this algebra with the involution
\begin{equation}\label{invf}
\left(t_{\{1,2,\ldots,n \}J}^{\wedge
n}\right)^*=(-1)^{\mathrm{card}(\{1,2,\ldots,n \}\cap
J)}(-q)^{l(J,J^c)}t_{\{n+1,n+2,\ldots,2n \}J^c}^{\wedge n},
\end{equation}
where $J^c=\{1,2,\ldots,2n \}\setminus J$ and
$l(J,J^c)=\mathrm{card}\{(j',j'')\in J \times J^c|\: j'>j''\}$. The
$*$-algebra arising this way is denoted by
$\mathrm{Pol}(\mathrm{Pl}_{n,2n})_q$. Let $t=t_{\{1,2,\ldots,n
\}\{n+1,n+2,\ldots,2n \}}^{\wedge n}$ and $x=tt^*$. Obviously, $x$
quasi-commutes with all the generators of
$\mathrm{Pol}(\mathrm{Pl}_{n,2n})_q$.

Let $\mathrm{Pol}(\mathrm{Pl}_{n,2n})_{q,x}$ be the localization of the
$*$-algebra $\mathrm{Pol}(\mathrm{Pl}_{n,2n})_q$ with respect to the
multiplicative system $x^\mathbb{N}$. The results of section 2 of
\cite{SSV2} can be used to prove the following statement intended to shed
some light to the way the $*$-algebra $\mathrm{Pol}(\mathrm{Pl}_{n,2n})_q$
works in our constructions.

\medskip

\begin{lemma}\label{emb} Let $J_{\alpha a}=\{a \}\cup
\{n+1,n+2,\ldots,2n\}\setminus \{2n+1-\alpha \}$. The map ${\mathcal
I}:z_a^\alpha \mapsto t^{-1}t^{\wedge n}_{\{1,2,\ldots,n \}J_a^\alpha}$
admits an extension up to an embedding of the $*$-algebras ${\mathcal
I}:\mathrm{Pol}(\mathrm{Mat}_n)_q \to
\mathrm{Pol}(\mathrm{Pl}_{n,2n})_{q,x}$
\end{lemma}

\medskip

The next step in proving theorem \ref{psi} is in producing an auxiliary
$*$-algebra $\mathscr{F}$, together with a homomorphism of $*$-algebras
$\varphi:\mathrm{Pol}(\mathrm{Pl}_{n,2n})_q \to \mathscr{F}$. Let
$\mathbb{C}[\mathrm{Pl}_{n,2n}]_q \subset
\mathrm{Pol}(\mathrm{Pl}_{n,2n})_q$ be the subalgebra generated by the
quantum minors $t^{\wedge n}_{\{1,2,\ldots,n \}J}$, $\mathrm{card}(J)=n$,
and $\mathbb{C}[\overline{\mathrm{Pl}_{n,2n}}]_q \subset
\mathrm{Pol}(\mathrm{Pl}_{n,2n})_q$ the subalgebra generated by the
quantum minors $t^{\wedge n}_{\{n+1,n+2,\ldots,2n \}J}$,
$\mathrm{card}(J)=n$.

$\mathscr{F}$ appears to be an extension of
$\mathbb{C}[\mathrm{Pl}_{n,2n}]_q$; it is given by adding the elements
$\eta$, $\eta^{-1}$ to the list of generators and the relations
$$
\eta \eta^{-1}=\eta^{-1}\eta=1,\qquad \eta t^{\wedge n}_{\{1,2,\ldots,n
\}J}=q^{-n}t^{\wedge n}_{\{1,2,\ldots,n \}J}\eta,\qquad \mathrm{card}(J)=n
$$
to the list of relations. (As it is well known, this list is exhausted by
the commutation relations and q-analogues of Pl\"ucker relations).

\medskip

\begin{lemma}
The map
\begin{equation}\label{phi1}
\varphi:t^{\wedge n}_{\{1,2,\ldots,n \}J}\mapsto t^{\wedge
n}_{\{1,2,\ldots,n \}J}
\end{equation}
\begin{equation}\label{phi2}
\varphi:t^{\wedge n}_{\{n+1,n+2,\ldots,2n \}J}\mapsto \eta t^{\wedge
n}_{\{1,2,\ldots,n \}J}
\end{equation}
admits an extension up to a homomorphism of algebras
$\varphi:\mathrm{Pol}(\mathrm{Pl}_{n,2n})_q \to \mathscr{F}$.
\end{lemma}

\smallskip

{\bf Proof.} The definitions imply the existence of a homomorphism
$\mathbb{C}[\mathrm{Pl}_{n,2n}]_q \to \mathscr{F}$ which satisfies
(\ref{phi1}) and a homomorphism
$\mathbb{C}[\overline{\mathrm{Pl}}_{n,2n}]_q \to \mathscr{F}$ which
satisfies (\ref{phi2}). The rest of the relations between the generators
of $\mathrm{Pol}(\mathrm{Pl}_{n,2n})_q$ are commutation relations. What
remains is to establish that the same commutation relations are also valid
for images of the elements $t^{\wedge n}_{\{1,2,\ldots,n \}J}$, $t^{\wedge
n}_{\{n+1,n+2,\ldots,2n \}I}$ with respect to $\varphi$. For that, we use
the R-matrix form of commutation relations between the matrix elements of
the corresponding fundamental representation of the Hopf algebra $U_q
\mathfrak{sl}_{2n}$.
$$
\mathrm{const}'(q,n)t^{\wedge n}_{\{n+1,n+2,\ldots,2n \}I}t^{\wedge
n}_{\{1,2,\ldots,n
\}J}=\sum_{\{I',J'|\mathrm{card}(I')=\mathrm{card(J')=n}\}}
R_{IJ}^{I'J'}t^{\wedge n}_{\{1,2,\ldots,n \}J'}t^{\wedge
n}_{\{n+1,n+2,\ldots,2n \}I'},
$$
$$
\mathrm{const}''(q,n)t^{\wedge n}_{\{1,2,\ldots,n \}I}t^{\wedge
n}_{\{1,2,\ldots,n
\}J}=\sum_{\{I',J'|\,\mathrm{card}(I')=\mathrm{card(J')=n}\}}
R_{IJ}^{I'J'}t^{\wedge n}_{\{1,2,\ldots,n \}J'}t^{\wedge n}_{\{1,2,\ldots,n
\}I'}.
$$

Here $(R_{IJ}^{I'J'})$ is the R-matrix from the right hand side of the well
known relation $RTT=TTR$ \cite{CP}, and the constants
$\mathrm{const}'(q,n)$, $\mathrm{const}''(q,n)$ describe the action of the
R-matrix in the left hand side of that relation. We are to prove now that
$\mathrm{const}''(q,n)=q^{-n}\mathrm{const}'(q,n)$. This follows from
$$
t_{JI}^{\wedge n}t_{IJ}^{\wedge n}=t_{IJ}^{\wedge n}t_{JI}^{\wedge n},\qquad
t_{II}^{\wedge n}t_{IJ}^{\wedge n}=q^nt_{IJ}^{\wedge n}t_{II}^{\wedge n},
$$
with $I=\{1,2,\ldots,n \}$, $J=\{n+1,n+2,\ldots,2n \}$.\hfill $\square$

\medskip

Equip $\mathscr{F}$ with an involution.

\medskip

\begin{lemma} i) There exists a unique involution $*$ in $\mathscr{F}$
such that
$$\eta^*=q^{-n(n-1)}\eta^{-1},$$
$$
\left(t^{\wedge n}_{\{1,2,\ldots,n
\}J}\right)^*=(-1)^{\mathrm{card}(\{1,2,\ldots,n \}\cap
J)}(-q)^{l(J,J^c)}\eta t^{\wedge n}_{\{1,2,\ldots,n \}J^c}.
$$

ii) The above homomorphism $\varphi:\mathrm{Pol}(\mathrm{Pl}_{n,2n})_q \to
\mathscr{F}$ is a homomorphism of $*$-algebras.
\end{lemma}

\smallskip

{\sc Remark.} The motives to deduce the latter equality are as follows:
\begin{multline*}
\left(t^{\wedge n}_{\{1,2,\ldots,n \}J}\right)^*=\left(\varphi
\left(t^{\wedge n}_{\{1,2,\ldots,n \}J}\right)\right)^*=\varphi \left(
\left(t^{\wedge n}_{\{1,2,\ldots,n \}J}\right)^*\right)=
\\ =(-1)^{\mathrm{card}(\{1,2,\ldots,n \}\cap
J)}(-q)^{l(J,J^c)}\varphi \left(t^{\wedge n}_{\{n+1,n+2,\ldots,2n
\}J^c}\right)= \\ =(-1)^{\mathrm{card}(\{1,2,\ldots,n \}\cap
J)}(-q)^{l(J,J^c)}\eta \cdot t^{\wedge n}_{\{1,2,\ldots,n \}J^c}.
\end{multline*}

\smallskip

{\bf Proof.} The uniqueness of the involution $*$ is straightforward, and
its existence follows from (\ref{invf}). More precisely, (\ref{invf}) and
the isomorphism $\mathbb{C}[\mathrm{Pl}_{n,2n}]\widetilde{\to}
\mathbb{C}[\overline{\mathrm{Pl}_{n,2n}}]$, $t^{\wedge n}_{\{1,2,\ldots,n
\}J} \mapsto t^{\wedge n}_{\{n+1,n+2,\ldots,2n \}J}$,
$\mathrm{card}(J)=n$, imply the existence of an antilinear
antiautomorphism $*:\mathscr{F} \to \mathscr{F}$ with the required
properties. The property $**=\mathrm{id}$ is to be verified separately:
$$
\left(t^{\wedge n}_{\{1,2,\ldots,n
\}J}\right)^{**}=q^{l(J,J^c)+l(J^c,J)}\eta t^{\wedge n}_{\{1,2,\ldots,n
\}J}\eta^*=q^{n^2-n}t^{\wedge n}_{\{1,2,\ldots,n \}J}\eta \eta^*=t^{\wedge
n}_{\{1,2,\ldots,n \}J}.
$$
Now the definition of the involution and (\ref{invf}) imply that $\varphi$
is a $*$-homomorphism. \hfill $\square$

\medskip

Note that $\varphi$ admits a unique extension up to a homomorphism
$\varphi_x:\mathrm{Pol}(\mathrm{Pl}_{n,2n})_{q,x}\to \mathscr{F}_x$ of the
localizations of $\mathrm{Pol}(\mathrm{Pl}_{n,2n})_q$ and $\mathscr{F}$
with respect to the multiplicative system $x^\mathbb{N}$, $x=tt^*$.

Here is the last element of our construction.

\bigskip

\begin{lemma} The map $z_a^\alpha \mapsto t^{-1}t^{\wedge
n}_{\{1,2,\ldots,n \}J_{a \alpha}}$, $a,\alpha=1,2,\ldots,n$, admits a
unique extension up to an embedding of $*$-algebras
$\mathrm{Pol}(S(\mathbb{U}))_q \to \mathscr{F}_x$.
\end{lemma}

\smallskip

{\bf Proof.} Note first that $l(J_{\alpha a},J_{\alpha a
}^c)=a+\alpha+(n-1)^2-2$. Hence
\begin{multline*}
\left(t^{-1}t^{\wedge n}_{\{1,2,\ldots,n \}J_{a \alpha}}\right)^*=
\\=-(-q)^{a+\alpha+(n-1)^2-2}\eta t^{\wedge
n}_{\{1,2,\ldots,n \}J_{\alpha a}^c}\left((-q)^{n^2}\eta t^{\wedge
n}_{\{1,2,\ldots,n \}\{1,2,\ldots,n \}}\right)^{-1}=
\\ =-(-q)^{a+\alpha-2n-1}{t^{\wedge n}_{\{1,2,\ldots,n \}J_{\alpha a}^c}}
\left(t^{\wedge n}_{\{1,2,\ldots,n \}\{1,2,\ldots,n \}}\right)^{-1}=
\\ =(-q)^{a+\alpha-2n}\left(t^{\wedge n}_{\{1,2,\ldots,n \}\{1,2,\ldots,n
\}}\right)^{-1}t^{\wedge n}_{\{1,2,\ldots,n \}J_{\alpha a}^c}=
\\ =(-q)^{a+\alpha-2n}\left(t^{-1}t^{\wedge n}_{\{1,2,\ldots,n
\}\{1,2,\ldots,n \}}\right)^{-1}\left(t^{-1}t^{\wedge n}_{\{1,2,\ldots,n
\}J_{\alpha a}^c}\right).
\end{multline*}

On the other hand, in the algebra $\mathbb{C}[\mathrm{Pl}_{n,2n}]_{q,t}$
defined as a localization of $\mathbb{C}[\mathrm{Pl}_{n,2n}]_q$ with respect
to the multiplicative system $t^\mathbb{N}$, one has the following
relations, established in \cite{SSV2}:
\begin{gather*}
t^{-1}t^{\wedge n}_{\{1,2,\ldots,n \}\{1,2,\ldots,n \}}=\det \nolimits_q
\mathbf{z},
\\ t^{-1}t^{\wedge n}_{\{1,2,\ldots,n \}J_{a \alpha}^c}=\det \nolimits_q
\mathbf{z}_a^\alpha.
\end{gather*}
What remains is to use the natural embedding
$\mathbb{C}[\mathrm{Pl}_{n,2n}]_{q,t}\to \mathscr{F}_x$. \hfill $\square$

\medskip

To complete the proof of theorem \ref{psi}, one has to observe that the
$*$-homomorphism $\mathrm{Pol}(\mathrm{Pl}_{n,2n})_{q,x}\to \mathscr{F}_x$
takes the images of $z_a^\alpha \in \mathrm{Pol}(\mathrm{Mat}_n)_q$ with
respect to the embedding $\mathrm{Pol}(\mathrm{Mat}_n)_q \hookrightarrow
\mathrm{Pol}(\mathrm{Pl}_{n,2n})_q$ to the images of the corresponding
elements $z_a^\alpha \in \mathrm{Pol}(S(\mathbb{U}))_q$ with respect to
the embedding $\mathrm{Pol}(S(\mathbb{U}))_q \hookrightarrow
\mathscr{F}_x$. Now theorem \ref{psi} is proved. \hfill $\square$

\medskip

\begin{proposition}
i) $\mathrm{Im}\psi=\mathrm{Pol}(S(\mathbb{U}))_q$.

ii) $\ker \psi \subset \mathrm{Pol}(\mathrm{Mat}_n)_q$ constitutes a $U_q
\mathfrak{sl}_{2n}$-submodule.
\end{proposition}

\smallskip

{\bf Proof.} The first statement follows from $(\det \nolimits_q
\mathbf{z})^{-1}=q^{n(n-1)}(\det \nolimits_q \mathbf{z})^*$, since
$z_a^\alpha \in \mathrm{Im}\psi$ for all $a,\alpha=1,2,\ldots,n$ and
$(\det \nolimits_q \mathbf{z})^*\in \mathrm{Im}\psi$. The second statement
involves the structures of $U_q \mathfrak{sl}_{2n}$-module algebras in
$\mathrm{Pol}(\mathrm{Mat}_n)_{q,x}$, $\mathscr{F}_x$, which are imposed
in an obvious way (the elements $\eta,\eta^{-1}\in \mathscr{F}_x$ are $U_q
\mathfrak{sl}_{2n}$-invariants). The statement to be proved follows now
from the fact that the homomorphisms $\mathrm{Pol}(\mathrm{Mat}_n)_q \to
\mathrm{Pol}(\mathrm{Pl}_{n,2n})_{q,x}$,
$\mathrm{Pol}(\mathrm{Pl}_{n,2n})_{q,x}\to \mathscr{F}_x$ are morphisms of
$U_q \mathfrak{sl}_{2n}$-modules. \hfill $\square$

\medskip

Recall the notation $E_i$, $F_i$, $K_i^{\pm 1}$ for the standard
generators of the Hopf algebra $U_q \mathfrak{sl}_{2n}$ and the notation
$U_q \mathfrak{su}_{n,n}$ for the Hopf $*$-algebra $(U_q
\mathfrak{sl}_{2n},*)$, with
\begin{align*}
E_n^*&=-K_nF_n,&F_n^*&=-E_nK_n^{-1},&(K_n^{\pm 1})^*&=K_n^{\pm 1},&
\\ E_j^*&=K_jF_j,&F_j^*&=E_jK_j^{-1},&(K_j^{\pm 1})^*&=K_j^{\pm1},
&\text{for}\quad j \ne n.
\end{align*}

Equip $\mathrm{Pol}(S(\mathbb{U}))_q$ with a structure of $U_q
\mathfrak{su}_{n,n}$-module algebra via the canonical isomorphism
$\mathrm{Pol}(S(\mathbb{U}))_q \simeq \mathrm{Pol}(\mathrm{Mat}_n)_q/\ker
\psi$. Thus, $S(\mathbb{U})_q$ is a homogeneous space of the quantum group
$SU_{n,n}$.

To conclude, introduce one more $U_q \mathfrak{su}_{n,n}$-module algebra.

It was noted earlier that $z_a^\alpha$, $a,\alpha=1,2,\ldots,n$, generate
the $*$-algebra $\mathrm{Pol}(S(\mathbb{U}))_q$. Consider its extension in
the class of $*$-algebras given by adding generator $t$, together with the
additional relations
$$
tt^*=t^*t,\quad tz_a^\alpha=q^{-1}z_a^\alpha t,\quad
t^*z_a^\alpha=q^{-1}z_a^\alpha t^*,\quad a,\alpha=1,2,\ldots,n.
$$
This algebra will be denoted by $\mathrm{Pol}(\widehat{S}(\mathbb{U}))_q$.
Our intention is to extend the structure of $U_q
\mathfrak{su}_{n,n}$-module algebra from $\mathrm{Pol}(S(\mathbb{U}))_q$
onto $\mathrm{Pol}(\widehat{S}(\mathbb{U}))_q$. This is accessible via
embedding the $*$-algebra $\mathrm{Pol}(\widehat{S}(\mathbb{U}))_q$ into
the $*$-algebra $\mathscr{F}_x$: $i:t \mapsto t$, $i:z_a^\alpha \mapsto
t^{-1}t^{\wedge n}_{\{1,2,\ldots,n \}J_{\alpha a}}$,
$a,\alpha=1,2,\ldots,n$. In fact, the $*$-subalgebra
$i(\mathrm{Pol}(\widehat{S}(\mathbb{U}))_q)$ contains all the elements
$t$, $t^{-1}t^{\wedge n}_{\{1,2,\ldots,n \}J}$, $\mathrm{card}(J)=n$, and
hence is a $U_q \mathfrak{su}_{n,n}$-module subalgebra. What remains is to
transfer this structure onto $\mathrm{Pol}(\widehat{S}(\mathbb{U}))_q$. It
follows, as, in \cite{SSV1}, that
\begin{gather*}
E_jt=F_jt=(K_j^{\pm 1}-1)t=0,\qquad j \ne n,
\\ F_nt=(K_n^{\pm 1}-1)t=0,\qquad E_nt=q^{-1/2}tz_n^n.
\end{gather*}

The $U_q \mathfrak{su}_{n,n}$-module algebra
$\mathrm{Pol}(\widehat{S}(\mathbb{U}))_q$ will be an essential tool in
producing the principal degenerate series of representations of the
quantum group $SU_{n,n}$.

\bigskip

\section{A $\mathbf{U_q \mathfrak{sl}_{2n}}$-invariant integral}

There exists a unique invariant integral on the quantum group $U_n$
normalized in such a way that $\displaystyle \int \limits_{(U_n)_q}\!\!\!
1d \nu=1$. Use the natural isomorphism of $*$-algebras
$\mathrm{Pol}(S(\mathbb{U}))_q \to \mathbb{C}[U_n]_q$, $z_a^\alpha \mapsto
q^{\alpha-n}z_a^\alpha$, $a,\alpha=1,2,\ldots,n$, to transfer this
integral onto $\mathrm{Pol}(S(\mathbb{U}))_q$. It is easy to demonstrate
that the linear functional we get this way
$$
\mathrm{Pol}(S(\mathbb{U}))_q \to \mathbb{C},\qquad f \mapsto \int
\limits_{S(\mathbb{U})_q}\!\!\! fd \nu,
$$
is not $U_q \mathfrak{sl}_{2n}$-invariant provided
$\mathrm{Pol}(S(\mathbb{U}))_q$ is equipped with the ordinary structure of
$U_q \mathfrak{sl}_{2n}$-module. This section is intended to prove the $U_q
\mathfrak{sl}_{2n}$-invariance of this integral with respect to a 'twisted'
structure of $U_q \mathfrak{sl}_{2n}$-module in
$\mathrm{Pol}(S(\mathbb{U}))_q$. It is custom to say 'these are not
functions that should be integrated, but highest order differential forms'.
The change of the structure of $U_q \mathfrak{sl}_{2n}$-module in
$\mathrm{Pol}(S(\mathbb{U}))_q$ is nothing more than a passage from
functions to highest order forms.

Recall the notation $E_i$, $F_i$, $K_i^{\pm 1}$ for the standard
generators of the Hopf algebra $U_q \mathfrak{sl}_{2n}$, and $U_q
\mathfrak{s}(\mathfrak{gl}_n \times \mathfrak{gl}_n)$ stands for its
subalgebra, generated by $K_n^{\pm 1}$, $E_j$, $F_j$, $K_j^{\pm 1}$, $j
\ne n$. It is easy to prove the $U_q \mathfrak{s}(\mathfrak{gl}_n \times
\mathfrak{gl}_n)$-invariance of the above integral $f \mapsto \displaystyle
\int \limits_{S(\mathbb{U})_q}\!\!\! fd \nu$.

We use the notation $U_q \mathfrak{s}(\mathfrak{u}_n \times
\mathfrak{u}_n)=(U_q \mathfrak{s}(\mathfrak{gl}_n \times
\mathfrak{gl}_n),*)$

\medskip

\begin{proposition}
There exists such a structure of $U_q \mathfrak{sl}_{2n}$-module in
$\mathrm{Pol}(S(\mathbb{U}))_q$ that

i) its restriction to the subalgebra $U_q \mathfrak{s}(\mathfrak{gl}_n
\times
\mathfrak{gl}_n)$ coincides with the restriction of the standard structure
of
$U_q \mathfrak{s}(\mathfrak{gl}_n \times \mathfrak{gl}_n)$-module in
$\mathrm{Pol}(S(\mathbb{U}))_q$,

ii) the integral $\mathrm{Pol}(S(\mathbb{U}))_q \to \mathbb{C}$, $f
\mapsto \displaystyle \int \limits_{S(\mathbb{U})_q}\!\!\! fd \nu$, is a
morphism of $U_q \mathfrak{sl}_{2n}$-modules.
\end{proposition}

\smallskip

{\bf Proof.} We first produce some structure of $U_q
\mathfrak{sl}_{2n}$-module in $\mathrm{Pol}(S(\mathbb{U}))_q$, and then
prove that it satisfies i) and ii).

Consider the $U_q \mathfrak{su}_{n,n}$-module $*$-algebra
$\mathrm{Pol}(\widehat{S}(\mathbb{U}))_q$ introduced in the previous
section. The construction implies that the $*$-algebra
$\mathrm{Pol}(\widehat{S}(\mathbb{U}))_q$ is generated by $t$, $t^*$, and
the elements of its $U_q \mathfrak{su}_{n,n}$-module subalgebra
$\mathrm{Pol}(S(\mathbb{U}))_q$. Consult \cite{SSV1} for explicit formulae
which describe the action of the standard generators $E_j$, $F_j$,
$K_j^{\pm 1}$, $j=1,2,\ldots,2n-1$; the action on the conjugate elements
is easily derivable from those formulae since
$$
(\xi f)^*=(S(\xi))^*f^*,\qquad \xi \in U_q \mathfrak{su}_{n,n},\qquad f \in
\mathrm{Pol}(\widehat{S}(\mathbb{U}))_q.
$$

The elements $t$, $t^*$, $x$ quasi-commute with all the generators $t$,
$z_a^\alpha$, $a,\alpha=1,2,\ldots,n$, of the $*$-algebra
$\mathrm{Pol}(\widehat{S}(\mathbb{U}))_q$. This allows one to consider the
localization $\mathrm{Pol}(\widehat{S}(\mathbb{U}))_{q,x}$ of this
$*$-algebra with respect to the multiplicative system $x^{\mathbb{N}}$ and
then to extend the structure of $U_q \mathfrak{su}_{n,n}$-module algebra
onto it. Equip $\mathrm{Pol}(\widehat{S}(\mathbb{U}))_{q,x}$ with a
bigrading:
$$
\deg t=(1,0),\qquad \deg t^*=(0,1),\qquad
\deg(z_a^\alpha)=\deg(z_a^\alpha)^*=(0,0),\qquad a,\alpha=1,2,\ldots,n.
$$
Obviously, the homogeneous components
$$
\mathrm{Pol}(\widehat{S}(\mathbb{U}))_{q,x}^{(i,j)}=\{f \in
\mathrm{Pol}(\widehat{S}(\mathbb{U}))_{q,x}|\:\deg f=(i,j)\}=t^{*i}\cdot
\mathrm{Pol}(S(\mathbb{U}))_q \cdot t^j
$$
form submodules of the $U_q \mathfrak{sl}_{2n}$-module
$\mathrm{Pol}(\widehat{S}(\mathbb{U}))_{q,x}$. Equip
$\mathrm{Pol}(S(\mathbb{U}))_q$ with a structure of $U_q
\mathfrak{su}_{n,n}$-module via the vector space isomorphism
$$
\mathrm{Pol}(S(\mathbb{U}))_q \to
\mathrm{Pol}(\widehat{S}(\mathbb{U}))_{q,x}^{(-n,-n)},\qquad f
\mapsto(t^*)^{-n}ft^{-n}.
$$
It follows from the $U_q \mathfrak{s}(\mathfrak{gl}_n \times
\mathfrak{gl}_n)$-invariance of $t$, $t^*$ that the new $U_q
\mathfrak{sl}_{2n}$-module structure in $\mathrm{Pol}(S(\mathbb{U}))_q$
coincides with the previous one on the subalgebra $U_q
\mathfrak{s}(\mathfrak{gl}_n \times \mathfrak{gl}_n)\subset U_q
\mathfrak{sl}_{2n}$. So, our integral is again $U_q
\mathfrak{s}(\mathfrak{gl}_n \times \mathfrak{gl}_n)$-invariant. What
remains is to prove that
\begin{align*}
\int \limits_{S(\mathbb{U})_q}\!\!\!(t^*)^nF_n((t^*)^{-n}ft^{-n})t^nd
\nu&=0,&
\\ \int \limits_{S(\mathbb{U})_q}\!\!\!(t^*)^nE_n((t^*)^{-n}ft^{-n})t^nd
\nu&=0,& f \in \mathrm{Pol}(S(\mathbb{U}))_q.
\end{align*}
Observe that the integral in question is a real linear functional
$\displaystyle \int \limits_{S(\mathbb{U})_q}\!\!\! f^*d
\nu=\overline{\displaystyle \int \limits_{S(\mathbb{U})_q}\!\!\! fd \nu}$,
and $\mathrm{Pol}(\widehat{S}(\mathbb{U}))_{q,x}$ is a $U_q
\mathfrak{su}_{n,n}$-module $*$-algebra. Thus, it suffices to prove the
following

\medskip

\begin{lemma}
For all $f \in \mathrm{Pol}(S(\mathbb{U}))_q \simeq
\mathrm{Pol}(\widehat{S}(\mathbb{U}))_{q,x}^{(0,0)}$ one has
\begin{equation}\label{inv1}
\int \limits_{S(\mathbb{U})_q}\!\!\!(t^*)^nF_n((t^*)^{-n}ft^{-n})t^nd
\nu=0.
\end{equation}
\end{lemma}

\smallskip

{\bf Proof.} Identify the subalgebra $\mathrm{Pol}(S(\mathbb{U}))_q$ with
its image under the embedding into $\mathscr{F}_x$ to note that $F_nt=F_n
\eta=0$, $t^*=(-q)^{n^2}\eta t \det_q \mathbf{z}$. Hence (\ref{inv1}) is
equivalent to
\begin{equation}\label{inv1a}
\int \limits_{S(\mathbb{U})_q}\!\!\!(\det \nolimits_q
\mathbf{z})^nF_n((\det \nolimits_q \mathbf{z})^{-n}f)d \nu=0.
\end{equation}
On the other hand, $K_n \det_q \mathbf{z}=q^2\det_q \mathbf{z}$, as one can
see from an explicit form for $K_nz_a^\alpha$, $a,\alpha=1,2,\ldots,n$.
Hence, by a virtue of the $U_q \mathfrak{s}(\mathfrak{u}_n \times
\mathfrak{u}_n)$-invariance of the integral in question (\ref{inv1a}) is
also equivalent to
\begin{equation}\label{inv1b}
\int \limits_{S(\mathbb{U})_q}\!\!\!(\det \nolimits_q \mathbf{z})^n
\widetilde{F}_n((\det \nolimits_q \mathbf{z})^{-n}f)d \nu=0,
\end{equation}
with $\widetilde{F}_n=K_nF_n$.

We start with proving (\ref{inv1a}), (\ref{inv1b}) in the special case
$f=z_n^n$. In view of the relation $F_nz_n^n=q^{1/2}$ one has
$$
(\det \nolimits_q \mathbf{z})^nF_n((\det \nolimits_q
\mathbf{z})^{-n}z_n^n)=(\det \nolimits_q \mathbf{z})^nF_n((\det \nolimits_q
\mathbf{z})^{-n}K_n^{-1}z_n^n)+q^{1/2}.
$$
Apply the explicit form for $F_nz_a^\alpha$, $K_n^{\pm 1}z_a^\alpha$,
$a,\alpha=1,2,\ldots,n$, found in \cite{SSV1} to establish that
$K_n^{-1}z_n^n=q^{-2}z_n^n$,
$$
(\det \nolimits_q \mathbf{z})^nF_n((\det \nolimits_q
\mathbf{z})^{-n})=q^{1/2}\frac{q^{2n}-1}{q^{-2}-1}(\det \nolimits_q
\mathbf{z})^{-1}\det \nolimits_q(\mathbf{z}_n^n).
$$
Hence
$$
(\det \nolimits_q \mathbf{z})^nF_n((\det \nolimits_q
\mathbf{z})^{-n}z_n^n)=q^{1/2}\frac{q^{2n}-1}{q^{-2}-1}(\det \nolimits_q
\mathbf{z})^{-1}\det \nolimits_q(\mathbf{z}_n^n)q^{-2}z_n^n+q^{1/2}.
$$
On the other hand, $(\det \nolimits_q \mathbf{z})^{-1}\det \nolimits_q
(\mathbf{z}_n^n)=(z_n^n)^*$. Thus to prove (\ref{inv1a}) in the special case
$f=z_n^n$, we need only to verify that
$$
\int \limits_{S(\mathbb{U})_q}\!\!\!(z_n^n)^*z_n^nd
\nu=\frac{1-q^2}{1-q^{2n}}.
$$
For that, it suffices to elaborate the following relations:
\begin{eqnarray*}
\int \limits_{S(\mathbb{U}_n)_q}\!\!\!(z_n^n)^*z_n^nd \nu&=&\int
\limits_{(U_n)_q}\!\!\!(z_n^n)^*z_n^nd \nu,
\\ \int \limits_{(U_n)_q}\!\!\!(z_n^n)^*z_n^nd
\nu&=&\frac{q^{-(n-1)}}{q^{n-1}+q^{n-3}+\ldots+q^{-(n-3)}+q^{-(n-1)}}=
\frac{1-q^2}{1-q^{2n}}.
\end{eqnarray*}
The first relation follows from the explicit form for the natural
isomorphism of $*$-algebras $\mathrm{Pol}(S(\mathbb{U}))_q \simeq
\mathbb{C}[U_n]_q$, and the second one is a consequence of the
orthogonality relations for the quantum group $SU_n$ (see \cite[p.
457]{CP}).

In the special case $f=z_n^n$, (\ref{inv1a}), and hence (\ref{inv1b}) are
proved.

Turn to the general case $f \in \mathrm{Pol}(S(\mathbb{U}))_q$.

Normally we work with admissible (see \cite{SSV1, SSV2}) modules over
quantum enveloping algebras. We also use the standard basis
$H_1,\ldots,H_{2n-1}$ of the Cartan subalgebra of $\mathfrak{sl}_{2n}$.

Let $H_0=\sum \limits_{j=1}^{n-1}jH_j+nH_n+\sum
\limits_{j=1}^{n-1}jH_{2n-j}$ be the element of the Cartan subalgebra of
the Lie algebra $\mathfrak{s}(\mathfrak{gl}_n \times
\mathfrak{gl}_n)\subset \mathfrak{sl}_{2n}$, which is in its center and is
normalized so that $[H_0,E_n]=2E_n$, $[H_0,F_n]=-2F_n$. Associate to every
admissible simple finite dimensional $U_q \mathfrak{s}(\mathfrak{gl}_n
\times \mathfrak{gl}_n)$-module $V$ a pair $(\boldsymbol{\lambda},k)$,
with $\boldsymbol{\lambda}$ being the highest weight of the corresponding
$U_q \mathfrak{sl}_n \otimes U_q \mathfrak{sl}_n$-module, and the number
$k$ is determined by $H_0v=2kv$, $v \in V$. It is easy to obtain the
following decomposition of the $U_q \mathfrak{s}(\mathfrak{gl}_n \times
\mathfrak{gl}_n)$-module $\mathrm{Pol}(S(\mathbb{U}))_q$ into a sum of
simple finite dimensional pairwise non-isomorphic $U_q
\mathfrak{s}(\mathfrak{gl}_n \times \mathfrak{gl}_n)$-modules:
$$
\mathrm{Pol}(S(\mathbb{U}))_q=
\bigoplus_{(\boldsymbol{\lambda},k)}V_{(\boldsymbol{\lambda},k)}.
$$
Assume that for some $f \in \mathrm{Pol}(S(\mathbb{U}))_q$ (\ref{inv1b})
fails. Our immediate intention is to prove that $f$ belongs to some
component $V_{(\boldsymbol{\lambda}',k')}$, and to find the corresponding
pair $(\boldsymbol{\lambda}',k')$.

The presence of the structure of $U_q \mathfrak{s}(\mathfrak{gl}_n \times
\mathfrak{gl}_n)$-module in $\mathrm{Pol}(S(\mathbb{U}))_q$ leads to a
structure of $U_q \mathfrak{s}(\mathfrak{gl}_n \times
\mathfrak{gl}_n)$-module in
$\mathrm{End}_{\mathbb{C}}\mathrm{Pol}(S(\mathbb{U}))_q$ and a morphism of
$U_q \mathfrak{s}(\mathfrak{gl}_n \times \mathfrak{gl}_n)$-module algebras
$U_q \mathfrak{sl}_{2n}\to
\mathrm{End}_{\mathbb{C}}\mathrm{Pol}(S(\mathbb{U}))_q$. Consider the
linear operator in $\mathrm{Pol}(S(\mathbb{U}))_q$:
$$
L:\psi \mapsto (\det \nolimits_q \mathbf{z})^n\widetilde{F}_n((\det
\nolimits_q \mathbf{z})^{-n}\psi).
$$
Prove that the $U_q \mathfrak{s}(\mathfrak{gl}_n \times
\mathfrak{gl}_n)$-submodule of the $U_q \mathfrak{s}(\mathfrak{gl}_n
\times \mathfrak{gl}_n)$-module
$\mathrm{End}_{\mathbb{C}}\mathrm{Pol}(S(\mathbb{U}))_q$ generated by $L$
is isomorphic to $V_{(\boldsymbol{\lambda}'',k'')}$, with
$\boldsymbol{\lambda}''=(0,\ldots,0,1,*,1,0,\ldots,0)$, $k''=-1$. In fact,
$\widetilde{F}_n$ generates a finite dimensional simple $U_q
\mathfrak{s}(\mathfrak{gl}_n \times \mathfrak{gl}_n)$-submodule of the
$U_q \mathfrak{s}(\mathfrak{gl}_n \times \mathfrak{gl}_n)$-module $U_q
\mathfrak{sl}_{2n}$ in view of the Serre relations for
$F_1,F_2,\ldots,F_{2n-1}$ and the relation
$$\mathrm{ad}_q(E_j)\widetilde{F}_n=0,\qquad j \ne n.$$
What remains is to apply the commutation relations between $K_j^{\pm 1}$,
$j=11,2,\ldots,2n-1$, and $\widetilde{F}_n$.

Now we are in a position to finish the proof that $f \in
V_{(\boldsymbol{\lambda}',k')}$ and to find $\boldsymbol{\lambda}'$, $k'$.
Use the natural embedding
$\mathrm{End}_{\mathbb{C}}\mathrm{Pol}(S(\mathbb{U}))_q \otimes
\mathrm{Pol}(S(\mathbb{U}))_q \to \mathrm{Pol}(S(\mathbb{U}))_q$ and the
fact that $L \in V_{(\boldsymbol{\lambda}'',k'')}$, $Lf \in \mathbb{C}1$,
to get a non-zero morphism of $U_q \mathfrak{s}(\mathfrak{gl}_n \times
\mathfrak{gl}_n)$-modules $V_{(\boldsymbol{\lambda}'',k'')}\otimes
\left(\bigoplus
\limits_{(\boldsymbol{\lambda},k)}V_{(\boldsymbol{\lambda},k)}\right)\to
\mathbb{C}$. This means that $k'=-k''$ and $\boldsymbol{\lambda}'$ is the
highest weight of the dual representation of $U_q \mathfrak{sl}_n \otimes
U_q \mathfrak{sl}_n$. Finally,
$\boldsymbol{\lambda}'=(1,0,\ldots,0,*,0,\ldots,0,1)$, $k'=1$.
Furthermore, $f$ is the lowest weight vector in
$V_{(\boldsymbol{\lambda}',k')}$ since $L$ is the highest weight vector in
$V_{(\boldsymbol{\lambda}'',k'')}$. It follows that $f=\mathrm{const}\cdot
z_n^n$, which contradicts to the relation (\ref{inv1b}) proved earlier for
$z_n^n$. \hfill $\square$

\bigskip

\section{On certain irreducible \boldmath $*$-representation of $U_q
\mathfrak{su}_{n,n}$}

Consider the embedding ${\cal I}:\mathbb{C}[\mathrm{Mat}_n]_q
\hookrightarrow \mathbb{C}[\mathrm{Pl}_{n,2n}]_{q,t}$ described in Lemma
\ref{emb} and another embedding
$$
{\cal J}:\mathbb{C}[\mathrm{Mat}_n]_q \to
\mathbb{C}[\mathrm{Pl}_{n,2n}]_{q,t},\qquad {\cal J}:f
\mapsto({\cal I}f)t^{-n}.
$$
It is easy to verify that ${\cal J}\mathbb{C}[\mathrm{Mat}_n]_q$ is a
submodule of the $U_q \mathfrak{sl}_{2n}$-module
$\mathbb{C}[\mathrm{Pl}_{n,2n}]_{q,t}$. Hence, there exists a unique
representation $\pi$ of $U_q \mathfrak{sl}_{2n}$ in the vector space
$\mathbb{C}[\mathrm{Mat}_n]_q$ such that
$$
\pi(\xi)f={\cal J}^{-1}\xi({\cal J}f),\qquad f \in
\mathbb{C}[\mathrm{Mat}_n]_q,\quad \xi \in U_q \mathfrak{sl}_{2n}.
$$
We refer to the previous sections for constructions of a $*$-homomorphism
$$
\psi:\mathrm{Pol}(\mathrm{Mat}_n)_q \to \mathrm{Pol}(S(\mathbb{U}))_q,\qquad
\psi:f \mapsto f|_{S(\mathbb{U})_q},
$$
and an invariant integral
$$
\mathrm{Pol}(S(\mathbb{U}))_q \to \mathbb{C},\qquad f \mapsto \int
\limits_{S(\mathbb{U})_q}fd \nu.
$$
These are to be used to equip the vector space
$\mathbb{C}[\mathrm{Mat}_n]_q$ with a Hermitian scalar product
$$
(f_1,f_2)\stackrel{\mathrm{def}}{=}\int
\limits_{S(\mathbb{U})_q}(f_2|_{S(\mathbb{U})_q})^*f_1|_{S(\mathbb{U})_q}d
\nu.
$$
The following lemma is a consequence of the definitions of $U_q
\mathfrak{sl}_{2n}$-module structures and invariance of the integral.

\medskip

\begin{lemma}\label{sr}
For all $\xi \in U_q \mathfrak{su}_{n,n}$, $f_1,f_2 \in
\mathbb{C}[\mathrm{Mat}_n]_q$ one has
$$(\pi(\xi)f_1,f_2)=(f_1,\pi(\xi^*)f_2).$$
\end{lemma}

\medskip

To rephrase, $\pi$ is a $*$-representation of $U_q \mathfrak{su}_{n,n}$ in
the pre-Hilbert space $\mathbb{C}[\mathrm{Mat}_n]_q$. The proof of
irreducibility for this representation uses the following

\begin{lemma}\label{ccl}
If for some $f \in \mathbb{C}[\mathrm{Mat}_n]_q$ and some
$m_1,m_2,\ldots,m_{2n-1}\in \mathbb{Z}$ one has
\begin{equation}\label{cc}
F_jf=0,\qquad K_j^{\pm 1}f=q^{\pm m_j}f,\qquad j=1,2,\ldots,2n-1,
\end{equation}
then $f \in \mathbb{C}1$.
\end{lemma}

\smallskip

{\bf Proof.} Due to $F_if=0$, $i \ne n$, $K_j^{\pm 1}f=q^{\pm m_j}f$,
$j=1,2,\ldots,2n-1$, it follows that, up to a constant complex multiple
$$
f=(z_n^n)^{k_n}\left(\det \nolimits_q \mathbf{z}_{\{n-1,n\}}^{\wedge 2
\{n-1,n\}}\right)^{k_{n-1}}\left(\det \nolimits_q
\mathbf{z}_{\{n-2,n-1,n\}}^{\wedge 3 \{n-2,n-1,n\}}\right)^{k_{n-2}}\cdots
(\det \nolimits_q \mathbf{z})^{k_1}
$$
for some $k_1,k_2,\ldots,k_n \in \mathbb{Z}_+$. Prove that $F_nf=0$
implies $k_1=k_2=\ldots=k_n=0$. Let $J$ be the two-sided ideal of
$\mathbb{C}[\mathrm{Mat}_n]_q$, determined by the 'non-diagonal'
generators $z_a^\alpha$, $a \ne \alpha$. Obviously, $K_n^{-1}J \subset J$,
$F_nJ \subset J$. This allows one to arrange computations modulo the ideal
$J$:
\begin{gather*}
F_nf=\mathrm{const}\prod_{a=1}^{n-1}(z_a^\alpha)^{\sum \limits_{i=1}^ak_i}
(z_n^n)^{\sum \limits_{i=1}^nk_i-1}\qquad (\mathrm{mod}\;J)
\\ \mathrm{const}=q^{1/2}\sum \limits_{j=1}^nq^{-2 \sum
\limits_{i=j+1}^nk_i}\cdot \frac{q^{-2k_j}-1}{q^{-2}-1}.
\end{gather*}
Now the isomorphism $\mathbb{C}[\mathrm{Mat}_n]_q/J \simeq
\mathbb{C}[z_1^1,z_2^2,\ldots,z_n^n]$ implies the relation
$\mathrm{const}=0$, that is, $k_1=k_2=\ldots=k_n=0$, which is just our
statement. \hfill $\square$

\medskip

\begin{proposition}\label{irr}
The representation $\pi$ is irreducible.
\end{proposition}

\smallskip

{\bf Proof.} Suppose $\pi$ is reducible. Then by lemma \ref{sr} it is a
sum of two non-trivial admissible subrepresentations in the subspaces
$L_1$, $L_2$: $\mathbb{C}[\mathrm{Mat}_n]_q=L_1 \oplus L_2$. Each of those
further decomposes as a sum of weight subspaces and, in particular,
possesses a lowest weight vector. This vector $f \ne 0$ satisfies the
relations
$$
\pi(F_j)f=0,\qquad \pi(K_j^{\pm 1})f=q^{\pm m_j}f,\qquad
j=1,2,\ldots,2n-1.
$$
On the other hand, in $\mathbb{C}[\mathrm{Pl}_{n,2n}]_{q,t}$ one has:
$$
\pi(K_j^{\pm 1})t=
\begin{cases}
q^{\mp 1}t,& j=n
\\ t,& j \ne n
\end{cases},\qquad
\pi(F_j)t=0,\qquad j=1,2,\ldots,2n-1.
$$
We conclude that $f$ is a solution of the equation system (\ref{cc}), so,
by a virtue of lemma \ref{ccl}, $f \in \mathbb{C}1$. Hence, $L_1 \supset
\mathbb{C}1$, $L_2 \supset \mathbb{C}1$, $L_1 \cap L_2 \ne 0$. This is a
contradiction coming from our assumption on reducibility of $\pi$. The
proposition is proved. \hfill $\square$

\bigskip

\section{The Cauchi-Szeg\"o integral representation}

Recall that $\mathbb{C}[\overline{\mathrm{Mat}}_n]_q \subset
\mathrm{Pol}(\mathrm{Mat}_n)_q$ stands for the subalgebra generated by
$(z_a^\alpha)^*$, $a,\alpha=1,2,\ldots,n$. We follow \cite{SSV2} in
introducing the algebra
$$
\mathbb{C}[\mathrm{Mat}_n \times \overline{\mathrm{Mat}}_n]_q=
\mathbb{C}[\mathrm{Mat}_n]_q^{\mathrm{op}}\otimes
\mathbb{C}[\overline{\mathrm{Mat}}_n]_q,
$$
with $\mathrm{op}$ indicating the change of a multiplication law to the
opposite one. This algebra is bigraded:
$$
\deg(z_a^\alpha \otimes 1)=(1,0),\qquad \deg(1 \otimes
(z_a^\alpha)^*)=(0,1).
$$
Its completion with respect to this bigrading is denoted by
$\mathbb{C}[[\mathrm{Mat}_n \times \overline{\mathrm{Mat}}_n]]_q$. The
elements of $\mathbb{C}[[\mathrm{Mat}_n \times
\overline{\mathrm{Mat}}_n]]_q$ are q-analogues for kernels of integral
operators, while the elements of the subalgebra $\mathbb{C}[\mathrm{Mat}_n
\times \overline{\mathrm{Mat}}_n]_q$ are q-analogues of polynomial
kernels.

We refer to \cite{SSV2} for a definition of the pairwise commuting 'kernels'
$$
\chi_k=\sum_{\genfrac{}{}{0mm}{}{J',J'' \subset \{1,2,\ldots,n\}}
{\mathrm{card}(J')=\mathrm{card}(J'')=k}}{z^{\wedge k}}_{J''}^{J'}\otimes
\left({z^{\wedge k}}_{J''}^{J'}\right)^*\in \mathbb{C}[\mathrm{Mat}_n \times
\overline{\mathrm{Mat}}_n]_q
$$
and a one-parameter family of the elements of $\mathbb{C}[[\mathrm{Mat}_n
\times \overline{\mathrm{Mat}}_n]]_q$ which includes the kernel
$$C_q=\prod_{j=0}^{n-1}\left(1+\sum_{k=1}^n(-q^{2j})^k \chi_k\right)^{-1}.$$
We write $C_q(\mathbf{z},\boldsymbol{\zeta}^*)$ instead of $C_q$,
$C_q(\mathbf{z},\boldsymbol{\zeta}^*)\cdot f(\boldsymbol{\zeta})$ instead
of $C_q \cdot(1 \otimes f)$, and $\displaystyle \int
\limits_{S(\mathbb{U})_q}f(\boldsymbol{\zeta})d \nu(\boldsymbol{\zeta})$
instead of $\displaystyle \int
\limits_{S(\mathbb{U})_q}f|_{S(\mathbb{U})_q}d \nu$, as it is custom in
the classical analysis. It is easy to demonstrate \cite[proposition
2.11]{SSV2} that in the formal passage to a limit $q \to 1$ one has
$C_q(\mathbf{z},\boldsymbol{\zeta}^*) \to
\det(1-\mathbf{z}\boldsymbol{\zeta}^*)^{-n}$.

We call $C_q \in \mathbb{C}[[\mathrm{Mat}_n \times
\overline{\mathrm{Mat}}_n]]_q$ the Cauchi-Szeg\"o kernel for the quantum
matrix ball.

The following result provides a motivation for this definition.

\medskip

\begin{theorem}\label{ir}
For any element $f \in \mathbb{C}[\mathrm{Mat}_n]_q$ one has
$$
f(\mathbf{z})=\int
\limits_{S(\mathbb{U})_q}C_q(\mathbf{z},\boldsymbol{\zeta}^*)
f(\boldsymbol{\zeta})d \nu(\boldsymbol{\zeta}).
$$
\end{theorem}

\smallskip

{\bf Proof.} Equip $\mathrm{Pol}(\mathrm{Mat}_n)_q$ with a grading:
$\deg(z_a^\alpha)=1$, $\deg(z_a^\alpha)^*=-1$, $a,\alpha=1,2,\ldots,n$. It
is easy to show that $\displaystyle \int
\limits_{S(\mathbb{U})_q}f|_{S(\mathbb{U})_q}d \nu=0$ for all $f$ with
$\deg(f)\ne 0$. Hence the integral operator
$$
T:\mathbb{C}[\mathrm{Mat}_n]_q \to \mathbb{C}[\mathrm{Mat}_n]_q,\qquad
T:f(z) \mapsto \displaystyle \int
\limits_{S(\mathbb{U})_q}C_q(\mathbf{z},\boldsymbol{\zeta}^*)
f(\boldsymbol{\zeta})d \nu(\boldsymbol{\zeta})
$$
is well defined. It follows from the invariance of the integral on the
Shilov boundary and the results of \cite[section 8]{SSV2} that
$\pi(\xi)T=T \pi(\xi)$ for all $\xi \in U_q \mathfrak{sl}_{2n}$.
Furthermore, $T1=1$. Hence $Tf=f$ for all $f \in
\mathbb{C}[\mathrm{Mat}_n]_q$ in view of proposition \ref{irr}. The
theorem is proved. \hfill $\square$

\medskip

Note that there exists another proof of theorem \ref{ir} which uses only
the orthogonality relations for the quantum group $U_n$ and one of the
Milne's relations for Schur's functions \cite{Mi}. We intend to publish
this proof in a future work.

\bigskip

\section{Appendix. Another description of the $*$-algebra
\boldmath$\mathrm{Pol}(S(\mathbb{U}))_q$ and its generalization onto the
case of rectangular matrices}

We produce a system of equations which distinguish the quantum Shilov
boundary from the quantum matrix space. This is certainly equivalent to
describing the $*$-algebra $\mathrm{Pol}(S(\mathbb{U}))_q$ in terms of
generators and relations.

\begin{proposition}\label{ng}
In $\mathrm{Pol}(S(\mathbb{U}))_q$ the following relations are valid:
\begin{equation}\label{zaa4}
\sum_{j=1}^nq^{2n-\alpha-\beta}z_j^\alpha(z_j^\beta)^*-\delta^{\alpha
\beta}=0,\qquad \alpha,\beta=1,2,\ldots,n.
\end{equation}
The left hand sides of these equations generate the kernel of the canonical
homomorphism $\psi:\mathrm{Pol}(\mathrm{Mat}_n)_q \to
\mathrm{Pol}(S(\mathbb{U}))_q$.
\end{proposition}

\smallskip

{\bf Proof.} The first statement is due to the fact that in
$\mathbb{C}[U_n]_q$
$$
\sum_{j=1}^nz_j^\alpha(z_j^\beta)^\star-\delta^{\alpha \beta}=0,\qquad
\alpha,\beta=1,2,\ldots,n,
$$
Let $\mathcal{I}$ be the two-sided ideal of $\mathrm{Pol}(\mathrm{Mat}_n)_q$
generated by the left hand sides of (\ref{ng}) and
$A=\mathrm{Pol}(\mathrm{Mat}_n)_q/\mathcal{I}$. One has to prove that the
canonical onto map $j:A \to \mathrm{Pol}(S(\mathbb{U}))_q$ is in fact an
isomorphism.

First prove that in $A$
\begin{equation}\label{dd*}
\det \nolimits_q \boldsymbol{z}\cdot(\det \nolimits_q
\boldsymbol{z})^*=q^{-n(n-1)}.
\end{equation}
Equip $A$ with a $\mathbb{Z}$-grading as follows:
$$
\deg(z_a^\alpha)=1,\qquad \deg((z_a^\alpha)^*)=-1,\qquad
a,\alpha=1,2,\ldots,n
$$
It was demonstrated in section 2 that the relation (\ref{dd*}) is valid in
$\mathrm{Pol}(S(\mathbb{U}))_q$. What remains is to use the fact that the
algebra $A$ is a $U_q \mathfrak{sl}_n$-module algebra, together with the
following statement:

\smallskip

\begin{lemma}\label{inv0}\hfill
\\ i) $\det \nolimits_q \boldsymbol{z}\cdot(\det \nolimits_q
\boldsymbol{z})^*\in A$ is a $U_q \mathfrak{sl}_n$-invariant of degree zero.
\\ ii) The subalgebra of $U_q \mathfrak{sl}_n$-invariants of degree zero in
$A$ is one-dimensional.
\end{lemma}

\smallskip

{\bf Proof} of lemma \ref{inv0}. The first statement is obvious. Turn to the
second statement. For any $U_q \mathfrak{sl}_n$-invariant of degree zero $f
\in A$ there exists such $U_q \mathfrak{sl}_n$-invariant of degree zero
$\widehat{f}\in \mathrm{Pol}(\mathrm{Mat}_n)_q$ that $f=j \cdot
\widehat{f}$.\footnote{Due to local finiteness of the $U_q
\mathfrak{sl}_n$-modules $\mathcal{I}$, $\mathrm{Pol}(\mathrm{Mat}_n)_q$,
$A$.} On the other hand, $\mathrm{Pol}(\mathrm{Mat}_n)_q$ is a free
$\mathbb{C}[\mathrm{Mat}_n \times \overline{\mathrm{Mat}}_n]$-module
generated by $1 \in \mathrm{Pol}(\mathrm{Mat}_n)_q$:
$$f_1 \otimes f_2: g_1 \otimes g_2 \mapsto g_1f_1 \otimes f_2g_2.$$
Consider $\widetilde{f}\in \mathbb{C}[\mathrm{Mat}_n \times
\overline{\mathrm{Mat}}_n]$ such that $\widehat{f}=\widetilde{f}\circ 1$. It
suffices to prove that $\widetilde{f}$ is in the subalgebra generated by
$\sum \limits_{j=1}^nq^{2n-\alpha-\beta}z_j^\alpha \otimes(z_j^\beta)^*$,
$\alpha,\beta=1,2,\ldots,n$. This is a consequence of the fact that
$\widetilde{f}$ is a $U_q \mathfrak{sl}_n$-invariant of degree zero, and a
q-analogue of the first main theorem of the theory of invariants \cite{GLR}.
\hfill $\square$

\medskip

It follows from (\ref{dd*}) that $\det_q \boldsymbol{z}\in A$ is invertible.
Hence one has a well defined homomorphism of algebras
$j':\mathrm{Pol}(S(\mathbb{U}))_q \to A$, $j':z_a^\alpha \mapsto
z_a^\alpha$, $a,\alpha=1,2,\ldots,n$. Obviously, $jj'=\mathrm{id}$. What
remains is to prove that $j'$ is onto, that is the elements $(\det_q
\boldsymbol{z})^{-1}$, $z_a^\alpha$, $a,\alpha=1,2,\ldots,n$, generate the
algebra $A$. It follows from (\ref{zaa4}) and well known properties of
quantum determinants that in $A$ one has
$$
(z_a^\alpha)^*=(-q)^{a+\alpha-2n}(\det \nolimits_q \boldsymbol{z})^{-1}\det
\nolimits_q \boldsymbol{z}_a^\alpha,\qquad a,\alpha=1,2,\ldots,n.
$$
So, $j$ is invertible. The proof of proposition \ref{ng} is complete. \hfill
$\square$

\medskip

Thus we get another description of the $*$-algebra
$\mathrm{Pol}(S(\mathbb{U}))_q$. It is to be used for producing a q-analogue
for the Shilov boundary of a unit ball in the space of rectangular matrices
$\mathrm{Mat}_{m,n}$, $m<n$.

The classical theory of Cartan domains \cite{Hua} provides a well known
procedure of producing the Shilov boundary $S(\mathbb{U}')$ of the unit ball
$\mathbb{U}'\subset \mathrm{Mat}_{m,n}$, together with the associated
Cauchi-Szeg\"o kernel for the unit ball $\mathbb{U}\subset \mathrm{Mat}_n$.
We are going to apply exactly this method in the quantum case.

Consider the $*$-subalgebra $\mathrm{Pol}(\mathrm{Mat}_{m,n})_q \subset
\mathrm{Pol}(\mathrm{Mat}_n)_q$ generated by $z_a^\alpha$, $\alpha>n-m$, and
the $*$-Hopf algebra $U_q \mathfrak{su}_{m,n}\subset U_q
\mathfrak{su}_{n,n}$ generated by $E_j$, $F_j$, $K_j^{\pm 1}$, $j<m+n$. It
is easy to demonstrate that $\mathrm{Pol}(\mathrm{Mat}_{m,n})_q$ is a $U_q
\mathfrak{su}_{m,n}$-module subalgebra of the $U_q
\mathfrak{su}_{m,n}$-module algebra $\mathrm{Pol}(\mathrm{Mat}_n)_q$. Up to
relabeling the generators, this $U_q \mathfrak{su}_{m,n}$-module algebra
coincides with the $*$-algebra of polynomials in the quantum matrix space
considered in \cite{SSV1, SSV2}.

Introduce the notation $\mathrm{Pol}(S(\mathbb{U}'))_q$ for the $*$-algebra
determined by its generators
$$z_a^\alpha,\qquad a=1,2,\ldots,n;\quad \alpha=n-m+1,n-m+2,\ldots,n,$$
the commutation relations loan from $\mathrm{Pol}(\mathrm{Mat}_{m,n})_q$
and the additional relations
$$
\sum_{j=1}^nq^{2n-\alpha-\beta}z_j^\alpha(z_j^\beta)^*=\delta^{\alpha,\beta},
\qquad \alpha,\beta=n-m+1,n-m+2,\ldots,n.
$$
Our construction implies a commutative diagram
$$
\begin{CD}
\mathrm{Pol}(\mathrm{Mat}_{m,n})_q @>>>\mathrm{Pol}(\mathrm{Mat}_n)_q \\
@VVV @VVV
\\ \mathrm{Pol}(S(\mathbb{U}'))_q @>>>\mathrm{Pol}(S(\mathbb{U}))_q.
\end{CD}
$$
where the vertical arrows stand for the homomorphisms of restriction of
'polynomials' onto the Shilov boundaries of the corresponding quantum balls.
The homomorphism $\mathrm{Pol}(S(\mathbb{U}'))_q \to
\mathrm{Pol}(S(\mathbb{U}))_q$ allows one to transfer an invariant integral
from $\mathrm{Pol}(S(\mathbb{U}))_q$ onto $\mathrm{Pol}(S(\mathbb{U}'))_q$.
The Cauchi-Szeg\"o kernel is defined as in section 5 except that in the
expression for $\chi_k$, the summing up in $J'\subset \{1,2,\ldots,n \}$ is
replaced with summing up in $J'\subset \{n-m+1,n-m+2,\ldots,n \}$. Another
result of section 5, a q-analogue of the integral Cauch-Szeg\"o
representation for the matrix ball $\mathbb{U}'$, implies a similar integral
representation for the quantum matrix ball in the space of rectangular
matrices.

\end{document}